\documentclass[12pt,a4paper]{article}
\usepackage{fancyhdr}
\headheight\baselineskip \headsep2\baselineskip \textwidth210mm
\advance\textwidth -2in \textheight287mm \advance\textheight-2in
\usepackage{epsfig}
\usepackage{graphicx}

\label{sec.sub.gdn}

\usepackage{amsfonts}
\usepackage{amsmath}

\begin{document}

\title{Adaptive two-stage sequential double sampling}

\author { Bardia Panahbehagh \footnote{Department of Mathematics, Kharazmi University, Tehran, Iran, Email
address: panahbehagh@khu.ac.ir}
 \ \ \ \  Afshin Parvardeh
\footnote{Department of Statistics, Isfahan University, Isfahan,
Iran}
 \ \ \ Babak Mohammadi \footnote{Research Center for Health, Aja University of Medical Science, Tehran,
 Iran}}

%}
\maketitle

\begin{abstract}
In many surveys inexpensive auxiliary variables are available that
can help us to make more precise estimation about the main
variable. Using auxiliary variable has been extended by regression
estimators for rare and cluster populations. In conventional
regression estimator it is assumed that the mean of auxiliary
variable in the population is known. In many surveys we don't have
such wide information about auxiliary variable. In this paper we
present a multi-phase variant of two-stage sequential sampling
based on an inexpensive auxiliary variable associated with the
survey variable in the form of double sampling. The auxiliary
variable will be used in both design and estimation stage. The
population mean is estimated by a modified regression-type
estimator with two different coefficient. Results will be
investigated using some simulations following Median and Thompson
(2004).
\end{abstract}

\section{Introduction}

Adaptive cluster sampling was introduced by Thompson (1990) as an
efficient sampling procedure for estimating totals and means of
rare and clustered populations. Because of lack of control on
final sample size and also problems that raise in performing the
design to define and use neighborhood, Salehi and Smith (2005)
proposed another adaptive design that does not require neighborhood and does not generate "edge units" in the sample, but
exploit clustering in the population to find rare
events with a reasonable bound for final sample size.\\
Panahbehagh et al. (2011) have investigated using auxiliary
variable just in the design in Adaptive two-stage sequential
sampling (ATS) in a real case study of fresh water mussel. Salehi
et al. (2013) with assuming that the population mean for auxiliary
variable is known, developed using auxiliary variable in
estimation stage by two modified regression estimator. Medina and
Thompson (2004) proposed a double sampling version of cluster
sampling named adaptive cluster double sampling for using
auxiliary information by regression estimator in Adaptive cluster
sampling. Here we are going to introduce a double sampling version
of adaptive two-stage sequential sampling for the situations which
there is no complete information about auxiliary variable. We present a
multi phase variant of adaptive two-stage sequential sampling that
obtained by combining the ideas of adaptive two-stage sequential
sampling and double
sampling.\\
In section 2 we introduce the design and respective notation.
Section 3 presents a regression type estimator with two different
coefficient with respective variance and variance estimator. In
section 4 we have some simulation to evaluate our design and in
section 5 the paper will be finished with some conclusion about
the design.

\section{Notation and sampling design}
Two-stage sequential sampling was initially proposed by Salehi and
Smith (2005) as a sample design for sampling rare and clustered
populations and then Brown et al. (2008) proposed an adaptive
version of that.
% Panahbehagh et al. (....) extended
%Two stage sequential sampling for multiple
%selection criteria and showed some asymptotic properties of sample mean in this design. \\
In adaptive two-stage sequential sampling (ATS) allocation of
second-stage effort among primary units is based on preliminary
information from the sampled primary units. Additional survey
effort is directed to primary units where the secondary units in
the initial sample have met a pre-specified criterion, or
condition (e.g., an individual from the rare population is
present). This design effectively over-samples primary units with
high values, compared with other primary units, a method
consistent with the approach recommended by Kalton and Anderson
(1986) for sampling rare populations.\\
Suppose we have a population of $N$ units partitioned into
$M$ primary sample units (PSU), each contain $N_h$ secondary
sample units (SSU). Let $\{(h,j), h=1, 2,..., M, j=1, 2, ...,
N_h\}$ denote the j-th unit in the h-th primary unit with an
associated measurement or count $y_{hj}$ and an auxiliary
variables $x_{hj}$. Then, $\bar{Y}_{N_{h}} =
\frac{1}{N_h}\sum_{j=1}^{N_h}y_{hj}$ is the mean of $y$ values
for h-th PSU and $\bar{Y}_{N} =
\frac{1}{N}\sum_{h=1}^{M}N_{h}\bar{Y}_{N_{h}}$ is the mean of the
whole population. $\bar{X}_{N_h}$ and $\bar{X}_{N}$ will define
the same.\\
The first stage of an adaptive two-stage sequential double
sampling design consists of selection a simple random sample $s$
of
size $m$ of $M$ PSUs.\\
The second stage contains two phases. The first phase consists of selecting an
initial conventional sample $s_{1h}$ of size $n_{1h}$ in h-th PSU
where $h\epsilon
s$.\\
Second phase consists of doing a sequential
sampling (like the second stage of a two-stage sequential
sampling) with a condition $C$, in each $s_{1h}, h\epsilon s$
based on auxiliary information or
both target and auxiliary information. The final sample in this phase named $s_{2h}$ with size $n_{2h}$\\
Then for each PSU we will have 3 estimators:
\begin{itemize}
\item $\bar{x}_{n_{1h}}$ that is an estimator for the inexpensive
variable $x$, and when $s_{h1}$ is gathered using SRSWOR we have
$\bar{x}_{n_{1h}}=\frac{1}{n_{1h}}\sum_{j\epsilon
s_{{1h}}}x_{hj}$. \item $\hat{t}_{yn_{2h}}$ and
$\hat{t}_{xn_{2h}}$ that are Murthy estimators for total of the
auxiliary and target variables in the population based on doing
ATS in h-th PSU, in the selected sample in the first phase
($s_{{1h}})$.
\end{itemize}

\section{A regression-type estimator with two different coefficient} The common estimator in
this design is Murthy estimator that is an unbiased estimator for
mean of the population. In this section we will introduce a
regression-type estimator for $\bar{Y}_{N}$ based on Murthy
estimator. Following Medina and Thompson (2004) "the regression
estimator will be constructed under the assumption that the
relationship between y and x can be modelled through a stochastic
regression model $\xi$ with mean
$E_\xi(y_{hj}|x_{hj})=x_{hj}\beta$ and variance
$var_{\xi}(y_{hj}|x_{hj})=\upsilon_{hj} \sigma^2$,
$\upsilon_{hj}=\varphi(x_{hj})$ where the function $\varphi$ is
assumed to be known. Throughout this paper we will consider the
role of regression model $\xi$ as the model-assisted survey
sampling approach (Sarndal et al., 1992); that is, we will suppose
that the relationship  between $y$ and $x$ is described reasonably
well by $\xi$, and consequently that the model can be used as an
instrument for constructing appropriate estimators of the
population parameters, but inference will not depend on the
assumed model and will rather be "design-based"".\\
Our main problem is estimating the $\bar{Y}_{N}$; however, because
of the regression model, it will also be required to
estimate the finite population regression parameter $\beta$.\\\\
Now we propose a known general form of regression estimator
(Sarendal et al., 1992, p.364) as below:
\begin{equation*}
\hat{\mu}_{reg}=\bar{y}_{n_{2}}+\beta(\bar{x}_{n_{1}}-\bar{x}_{n_{2}})
\end{equation*}
where
\begin{equation*}
\bar{y}_{n_{2}}=\frac{1}{N}\sum_{h \epsilon s}
a_h\frac{\hat{t}_{yn_{2h}}}{\pi_h}, a_h=\frac{N_h}{n_{1h}}
\end{equation*}
that $\hat{t}_{yn_{2h}}$ is Murty estimator in $s_{1h}$ and
$\pi_h$ is probability of choosing h-th PSU in the first stage of
sampling. $\bar{x}_{n_{2}}$ is defined the same but for $x$. Also
\begin{equation*}
\bar{x}_{n_{1}}=\frac{1}{N}\sum_{h \epsilon
s}a_h\frac{t_{xn_{1h}}}{\pi_h}, t_{xn_{1h}}=\sum_{j \epsilon
s_{1h}} x_{hj}
\end{equation*}
 $\beta$ is a parameter and if it is unknown we should first
estimate it. Two reasonable candidates for $\beta$ are (see Salehi
et al. 2013)
\begin{equation*}
\hat{\beta}_1=\frac{\hat{t}_{{xy}{n_{2}}}-N\bar{y}_{n_{2}}\bar{x}_{{n_{2}}}}{
\hat{t}_{x^2}{_{n_{2}}}-N\bar{x}^2_{n_{2}}}
\end{equation*}
 and
\begin{eqnarray*}
\hat{\beta}_{o}=\frac{\hat{cov}(\bar{y}_{n_{2}},\bar{x}_{n_{2}})}{\hat{var}(\bar{x}_{n_{2}})}
\end{eqnarray*}
that are estimators of the conventional and the optimal regression
coefficient in ATS as below
\begin{eqnarray*}
\beta_1&=&\frac{\sum_{j=1}^{N}y_jx_j-N\bar{X}_{N}
\bar{Y}_{N}}{\sum_{j=1}^{N}x^2_j-N\bar{X}_{N}^2},\\
\beta_{o}&=&\frac{{cov}(\bar{y}_{n_{2}},\bar{x}_{n_{2}})}{{var}(\bar{x}_{n_{2}})}
\end{eqnarray*}

where $\hat{t}_{{xy}{n_{2}}}$ and $\hat{t}_{{x^2}{n_{2}}}$ are
unbiased Murthy estimators of the total of $xy$ and $x^2$ in the
population
respectively based on the design.\\

\subsection{Expectation and Variance of the estimators}
 With assuming $\hat{\beta}\simeq \beta$ according
 to the stages and the phases of the design we have
\begin{eqnarray*}
E(\hat{\mu}_{reg})&=&E_1E_2E_3(\hat{\mu}_{reg})
\end{eqnarray*}
and
\begin{eqnarray*}
var(\hat{\mu}_{reg})&=&V_1E_2E_3(\hat{\mu}_{reg})+E_1V_2E_3(\hat{\mu}_{reg})+E_1E_2V_3(\hat{\mu}_{reg})\\
&=&part1 + part2 + part3
\end{eqnarray*}
where $E$ and $V$ denote expectation and variance and the indexes
$1,2,3$, consist of first stage ($s$), second stage ($s_{1h}$) and
adaptive sampling in second stage ($s_{2h}$) respectively. Then
with SRSWOR in the first and second stage we have (see appendix
A):
\begin{eqnarray*}
E(\hat{\mu}_{reg})=\bar{Y}_{N}
\end{eqnarray*}
and
\begin{eqnarray*}
var(\hat{\mu}_{reg})=\frac{1}{N^2}[M^2(1-\frac{m}{M})\frac{S^{2}_{ty_N}}{m}+\frac{M}{m}\sum_{h=1}^{M}
N^2_{h}(1-\frac{n_{1h}}{N_h})\frac{S^2_{y_{N_h}}}{n_{1h}}\\+\frac{M}{m}\sum_{h=1}^{M}a^2_hE_2(V_3(\bar{y}_{n_{2h}})+\beta^2V_3(\bar{x}_{n_{2h}})-2\beta
C_3(\bar{y}_{n_{2h}},\bar{x}_{n_{2h}}))]
\end{eqnarray*}
where
\begin{eqnarray*}
S^{2}_{ty_N}=\frac{1}{M-1}\sum_{h=1}^{M}(t_{y_{N_h}}-\bar{t}_y)^2,
\bar{t}_y=\frac{1}{M}\sum_{h=1}^{M}t_{y_{N_h}}
\end{eqnarray*}
and
\begin{eqnarray*}
S^{2}_{y_{N_h}}=\frac{1}{N_h-1}\sum_{j=1}^{N_h}(y_{hj}-\bar{Y}_{N_h})^2,
\bar{Y}_{N_h}=\frac{1}{N_h}\sum_{j=1}^{N_h}y_{hj}
\end{eqnarray*}

\subsection{An estimator for $\beta_o$}
To calculate $\beta_o$, we have (see Appendix A)
\begin{eqnarray*}
var(\bar{x}_{n_{2}})=\frac{1}{N^2}[M^2(1-\frac{m}{M})\frac{S^{2}_{tx_N}}{m}+\frac{M}{m}\sum_{h=1}^{M}
N^2_{h}(1-\frac{n_{1h}}{N_h})\frac{S^2_{x_{N_h}}}{n_{1h}}\\+\frac{M}{m}\sum_{h=1}^{M}a^2_hE_2V_3(\hat{t}_{xn_{2h}})]
\end{eqnarray*}
and

\begin{eqnarray*}
cov(\bar{x}_{n_{2}},\bar{y}_{n_{2}})=\frac{1}{N^2}[M^2(1-\frac{m}{M})\frac{S_{tx_Nty_N}}{m}+\frac{M}{m}\sum_{h=1}^{M}
N^2_{h}(1-\frac{n_{1h}}{N_h})\frac{S_{x_{N_h}y_{N_h}}}{n_{1h}}\\+\frac{M}{m}\sum_{h=1}^{M}a^2_hE_2C_3(\hat{t}_{xn_{2h}},\hat{t}_{yn_{2h}})]
\end{eqnarray*}

But $\beta_o$ is a parameter yet and should be estimated by sample
information. For estimating variance and covariance terms, since
(see Appendix B)

\begin{eqnarray*}
E(\frac{1}{N^2}M^2(1-\frac{m}{M})\frac{\hat{S}^{2}_{tx_N}}{m})=\frac{1}{N^2}[M^2(1-\frac{m}{M})\frac{S^{2}_{tx_N}}{m}+\frac{M}{m}\sum_{h=1}^{M}
N^2_{h}(1-\frac{n_{1h}}{N_h})\frac{S^2_{x_{N_h}}}{n_{1h}}\\+\frac{M}{m}\sum_{h=1}^{M}a^2_hE_2V_3(\hat{t}_{xn_{2h}})-\sum_{h=1}^{M}
N^2_{h}(1-\frac{n_{1h}}{N_h})\frac{S^2_{x_{N_h}}}{n_{1h}}-\sum_{h=1}^{M}a^2_hE_2V_3(\hat{t}_{xn_{2h}})]
\end{eqnarray*}

and (see Appendix B)

\begin{eqnarray*}
E(\hat{S}^2_{x_{N_h}})=S^2_{x_{N_h}}-\frac{1}{n_{1h}(n_{1h}-1)}E_2V_3(\hat{t}_{xn_{2h}})
\end{eqnarray*}
where
\begin{eqnarray*}
\hat{S}^2_{x_{N_h}}=\frac{1}{n_{1h}-1}[\hat{t}_{x^2_{n_{2h}}}-\frac{\hat{t}^2_{x_{n_{2h}}}}{n_{1h}}],
\end{eqnarray*}
an reasonable estimator for $var(\bar{x}_{n_{2}})$ would be
\begin{eqnarray*}
\hat{var}(\bar{x}_{n_{2}})=\frac{1}{N^2}[M^2(1-\frac{m}{M})\frac{\hat{S}^{2}_{\hat{t}x_N}}{m}+\frac{M}{m}\sum_{h\in
s}
N^2_{h}(1-\frac{n_{1h}}{N_h})\frac{\hat{S}^2_{x_{N_h}}}{n_{1h}}\\+\frac{M}{m}\sum_{h\in
s}a^2_h\frac{n_{1h}(N_h-1)}{N_h(n_{1h}-1)}\hat{V}_3(\hat{t}_{xn_{2h}})]
\end{eqnarray*}

and for $cov(\bar{x}_{n_{2}},\bar{y}_{n_{2}})$ in a similar way we
have

\begin{eqnarray*}
\hat{cov}(\bar{x}_{n_{2}},\bar{y}_{n_{2}})=\frac{1}{N^2}[M^2(1-\frac{m}{M})\frac{\hat{S}_{\hat{t}x_N\hat{t}y_N}}{m}+\frac{M}{m}\sum_{h\in
s}
N^2_{h}(1-\frac{n_{1h}}{N_h})\frac{\hat{S}_{x_{N_h}y_{N_h}}}{n_{1h}}\\+\frac{M}{m}\sum_{h\in
s}a^2_h\frac{n_{1h}(N_h-1)}{N_h(n_{1h}-1)}\hat{C}_3(\hat{t}_{xn_{2h}},\hat{t}_{yn_{2h}})]
\end{eqnarray*}
and then an reasonable and asymptotically unbiased estimator for
$var(\hat{\mu}_{reg})$ is (see appendix B)
\begin{eqnarray*}
\hat{var}(\hat{\mu}_{reg})=\frac{1}{N^2}[M^2(1-\frac{m}{M})\frac{\hat{S}^{2}_{ty_N}}{m}+\frac{M}{m}\sum_{h\epsilon
s}
N_{h}(1-\frac{n_{1h}}{N_h})\frac{\hat{S}^2_{y_{N_h}}}{n_{1h}}\\+\frac{M}{m}\sum_{h
\epsilon
s}a^2_h\frac{n_{1h}(N_h-1)}{N_h(n_{1h}-1)}\hat{V}_3(\hat{t}_{yn_{2h}})+\frac{M^2}{m^2}\sum_{h
\epsilon
s}a^2_h(\hat{\beta}^2\hat{V}_3(\bar{x}_{n_{2h}})-2\hat{\beta}\hat{C}_3(\hat{t}_{yn_{2h}},\hat{t}_{xn_{2h}}))]
\end{eqnarray*}
where
\begin{eqnarray*}
\hat{S}^2_{ty_{N}}=\frac{1}{m-1}\sum_{h \epsilon s}
(a_h\hat{t}_{y{n_{2h}}}-\bar{\hat{t}}_{y{n_2}})^2,
\bar{\hat{t}}_{y{n_2}}=\frac{1}{m}\sum_{h \epsilon s}
a_h\hat{t}_{y{n_{2h}}},
\end{eqnarray*}

\begin{eqnarray*}
\hat{S}^2_{y_{N_h}}=\frac{1}{n_{1h}-1}(\hat{t}_{y^2n_{2h}}-\frac{\hat{t}_{yn_{2h}}^2}{n_{1h}}),
\end{eqnarray*}

\begin{eqnarray*}
\hat{V}_3(\hat{t}_{yn_{2h}})=n_{1h}\{p_{2h1}[\frac{(n_{1h}-1)(l_{2h1}-1)}{n_{2h1}-1}+\frac{l_{2h}-1}{l_{2h}}
((1-p_{2h1})\frac{n_{1h}-n_{2h1}}{n_{2h1}-1}-n_{1h}p_{2h1})]S^2_{x,2hc}\\
+p_{2h1}(1-p_{2h1})\frac{n_{1h}-n_{2h1}}{n_{2h1}-1}
(\overline{x}_{s_{2hc}}-\overline{x}_{s_{2hc'}})^2\\
+(1-p_{2h1})[\frac{(n_{1h}-1)(n_{2h1}-l_{2h1}-1)}{n_{2h1}-1}+\frac{n_{2h}-l_{2h}-1}{n_{2h}-l_{2h}}
(p_{2h1}\frac{n_{1h}-n_{2h1}}{n_{2h1}-1}-n_{1h}(1-p_{2h1}))]S^2_{x,2hc'}\}
\end{eqnarray*}
and
\begin{eqnarray*}
\hat{C}_3(\hat{t}_{yn_{2h}},\hat{t}_{xn_{2h}})=n_{1h}\{p_{2h1}[\frac{(n_{1h}-1)(l_{2h1}-1)}{n_{2h1}-1}+\frac{l_{2h}-1}{l_{2h}}
((1-p_{2h1})\frac{n_{1h}-n_{2h1}}{n_{2h1}-1}-n_{1h}p_{2h1})]S_{xy,2hc},\\
+p_{2h1}(1-p_{2h1})\frac{n_{1h}-n_{2h1}}{n_{2h1}-1}
(\overline{y}_{s_{2hc}}-\overline{y}_{s_{2hc'}})(\overline{x}_{s_{2hc}}-\overline{x}_{s_{2hc'}})\\
+(1-p_{2h1})[\frac{(n_{1h}-1)(n_{2h1}-l_{2h1}-1)}{n_{2h1}-1}+\frac{n_{2h}-l_{2h}-1}{n_{2h}-l_{2h}}
(p_{2h1}\frac{n_{1h}-n_{2h1}}{n_{2h1}-1}-n_{1h}(1-p_{2h1}))]S_{xy,2hc'}\}
\end{eqnarray*}
where $l_{2h}$ and $l_{2h1}$ are number of SSU in h-th PSU in the
total sample and the primary sample that satisfy in condition $C$
respectively and $n_{2h1}$ is the size of the primary sample in
h-th PSU to perform an ATS. \\ Also $\hat{\beta}$ can be even
$\hat{\beta}_o$ or $\hat{\beta}_1$ with respect to the
coefficient that is used in the regression estimator.\\
\section{Monte Carlo study}
In this section, following Medina and Thompson (2004), to
investigate the design and the estimators, we simulated two
populations with different features, and each of them with
different auxiliary variables. Each population obtained by dividing a
$20\times 20$ unit square into N=400 unit quadrates, partitioned
in 4 PSUs with equal size. We associated with the unit or quadrat (SSU) $u_{hj}$ a
vector $(y_{hj}, x_{hj}, z_{hj})$, where $y_{hj}$, $x_{hj}$ and
$z_{hj}$ denote the j-th value of the survey variable y in h-th PSU,
and the values of two auxiliary variables x and z. Information
about populations are in table 1. We generated the spatial pattern
following the Poisson cluster process. The number of clusters was
selected from a Poisson distribution, and cluster centers were
randomly located throughout the site. Individuals within the
cluster were located around the cluster center at a random
distance following an exponential distribution and a random
direction following a uniform distribution. Also we used another
variable in simulations, say $w$, where $w$ was a binary variable
defined as $w_{hj}=1$ if $y_{hj}>0$ and $w_{hj}=0$ otherwise.\\

\subsection{Expectation of the designs costs}
To compare fairly the design, we have derived analytic formula for
the expectation of adaptive two-stage sequential double sampling
cost and its conventional sampling counterpart with equal effort.
 The sampling designs considered in this study were compared using the
expected value of the cost function,
$Cost_T=c_{aux}n_{aux}+c_{tar}n_y$, where $Cost_T$ is the total
cost, $c_{aux}$ and $c_{tar}$ were the per element costs of
measuring the auxiliary variable and the target variable,
respectively, and $n_{aux}$ and $n_y$ were the total numbers of
measurements of the
auxiliary variable and the target variable.\\
In $Cost_T$ formula, $c_{aux}$, $c_{tar}$ and $n_{aux}$ are
constant and just $n_y$ is variable. Then we have

\begin{eqnarray*}
E(Cost_T)=c_{aux}n_{aux}+c_{tar}E(n_y)
\end{eqnarray*}

 Let $L_{h}$, $l_{1h}$ and
$l_{2h}$ be the number of units satisfying condition $C$ in the
h-th PSU, in the first phase and the second phase of the second
stage of sampling respectively.
 Furthermore, let $L^{(r)}_h$, $l^{(r)}_{1h}$ and $l^{(r)}_{2h}$ are the
number of  rare units in the h-th PSU, in the first phase and the
second phase of the second stage of sampling, respectively. Then
\begin{eqnarray*}
n_y=\sum_{h=1}^M (n_1+dl_{1h})I_h
\end{eqnarray*}

where  $n_1$ is the  size of initial sample in ATSD in the h-th PSU
and $I_{h}$ is an indicator function that takes ``1" when the h-th
PSU is selected in the first stage and ``0" otherwise. Since
 $l_{1h}|_{I_h=1,n_{1h},L_{1h}}\sim HG(n_1,n_{1h},L_{1h})$
 where HG denote Hypergeometric distribution and $L_{1h}$ denotes number of unites that satisfy in
 condition $C$ in the selected sample of size $n_{1h}$, we have (with
 $p^c_h=\frac{L_h}{N_h}$)
\begin{eqnarray*}
E(l_{1h})=E(E(l_{1h}|I_h=1,n_{1h},L_{1h}))=E(n_1\frac{L_{1h}}{n_{1h}})
=n_1\frac{E(L_{1h})}{n_{1h}}=n_1\frac{n_{1h}p^c_h}{n_{1h}}=n_1p^c_h
\end{eqnarray*}
and
\begin{eqnarray*}
E(I_h)=\frac{m}{M},
\end{eqnarray*}
therefore
\begin{eqnarray*}
E(n_y)=\frac{m}{M}\sum_{h=1}^M (n_1+dn_1p^c_h)
\end{eqnarray*}

Then to have a fair comparison,
\begin{itemize}
\item when we want to execute a ATS
with just target variable with equal effort, we should set
$E(Cost_T)=E(Cost_{ATS})=n_{y_{ATS}}c_{tar}$, therefore we should
set the initial sample size in each PSU as
\begin{eqnarray*}
n_1\simeq \frac{E(Cost_T)}{c_{tar}[\frac{m}{M}\sum_{h=1}^M
(1+d_{ATS}p_h)]}
\end{eqnarray*}
 where $p_h$ is percentage of
rare units in h-th PSU
\item for two-stage sampling we set
\begin{eqnarray*}
n\simeq \frac{E(Cost_T)}{mc_{tar}}
\end{eqnarray*}
\item for SRSWOR we set
\begin{eqnarray*}
n\simeq \frac{E(Cost_T)}{c_{tar}}
\end{eqnarray*}
\item and for regression in Two-stage Double sampling because this
design use as much as ATSD of auxiliary variables (i.e. $mn_{h1}$)
then it is enough to set number of target sample size that should
be taken in each selected PSU as
\begin{eqnarray*}
n_{y_{tR}}\simeq \frac{E(n_y)}{m}.
\end{eqnarray*}
(note that the symbols $n_.$, $c_.$ and $Cont_.$, define the same
as $ATSD$, but $"."$ is replaced with a proper symbol according to
respective design).\\

\end{itemize}

We used 4 designs and 7 estimators in the simulations:
\begin{itemize}
\item Adaptive two-stage sequential double sampling that used both
target and auxiliary variables for condition $C$ and the
estimators in this design were
\begin{itemize}
\item Two regression estimators with $\beta_o$ and $\hat{\beta}_o$
named $RegO$ and $Regopt$ \item Two regression estimator with
$\beta_1$ and $\hat{\beta}_1$ named $Reg1$ and $Regb1$
\end{itemize} \item Adaptive two-stage sequential sampling that used just target variables as
 condition $C$ with Murthy estimator named $ATS$

\item Double sampling that used simple random sampling for both
two phases with regression estimator that used sample mean named
$Regs$
 \item Simple random sampling without replacement named $\bar{y}_s$
\end{itemize}
Efficiency was defined as
$eff(\hat{y}_u)=\frac{var(\bar{y}_s)}{MSE(\hat{y}_u)}$ and
relative bias was defined as
$rbias(\hat{y}_u)=\frac{E(\hat{y}_u)-\bar{Y}_{N}}{\bar{Y}_{N}}$.
Condition $C$ was defined as " the respective SSU is nonempty" and
it depends on the design that used just target or both target and
auxiliary variables. Also in the iterations when it was not
possible to calculate the respective $\hat{\beta}$, we used
$\bar{y}_{n_{2}}$ for respective regression estimator.
\\
Two values for the ratio of costs $c_{aux}/c_{tar}$ were
considered, $c_{tar}/c_{aux}=5$ and $c_{tar}/c_{aux}=10$. In each
case the parameters were chosen such that total cost for
all the designs be almost the same.\\

\begin{table}[h]
\begin{center}
  \caption{ the feature of the populations.\label{tab}}
  \end{center}
\begin{center}
\begin{tabular}{c|ccc|ccc}
&&population1&&&population2&\\
&&rare and cluster&&&not so rare but cluster&\\ \hline &y&x&z&y&x&z\\
\hline
mean&0.19&0.40&0.42&0.50&0.47&0.46\\
variance&1.68&3.14&3.35&4.85&3.52&2.99\\
correlation with y&&0.89&0.52&&0.93&0.51\\
\end{tabular}%
\end{center}

\end{table}
The results for Population1 are in table 2 and 3 and can be
summarized as follow. For efficiency in the case of
$c_{tar}/c_{aux}=10$, adaptive two-stage sequential double
sampling (ATSD) is appropriate (albeit $w$ shows no regular pattern). In
the case of $c_{tar}/c_{aux}=5$, just for enough high correlation
(using $x$) we can trust ATSD and for low correlations ordinary
ATS (that expense all the costs for sampling target variable using ATS)
is more appropriate than others. Gain in efficiency for $Regb$ and $Regopt$ relative to $Regs$ is considerable .\\
Also results show that $n_{1h}$ and portion of
$\frac{n_{2h1}}{n_{1h}}$ are two important factors to improve the
efficiency of the estimators in the design. It is expected that
with increasing $n_{1h}$ the efficiency increases, but it is
interesting that bigness of $\frac{n_{2h1}}{n_{1h}}$ can amend
smallness of $n_{1h}$. For comparing $Regopt$ and $Regb1$
according to the results for high correlation $Regopt$ has better
performance than $Regb1$ and when the correlations are low we can
trust $Regb1$ more than $Regopt$. It could be a result of complexity of the $Regopt$ formula (see Salehi et. al. 2013)\\
For unbiasedness the results
can be summarized as follows. $n_{1h}$ is one of the important
factor that with increasing it, the bias decreases and the next
important factors are $m$ and $n_{2h1}$. In the cases that $ATSD$
is better than $ATS$, $Regopt$ has better (or at least equivalent)
performance than $Regb1$ and in some cases (for example the cases
that the correlations are low and $Regb1$ has good performance in
efficiency) bias of $Regopt$ is substantially smaller than $Regb1$
and bias of $Regb1$ is almost unacceptable. Also the amount of
bias for both $Regopt$ and $Regb1$ are unacceptable when our
auxiliary
variable is $w$ with $n_{1h}=50$.\\
Then for population1 with looking at efficiency and unbiasedness
together, in the cases that $ATSD$ is better than $ATS$, we can
trust $Regopt$ more than $Regb1$ and $Regopt$ can be our first
candidate to estimate the parameters.

\begin{table}[h]
\begin{center}
  \caption{ efficiency and relative bias of the estimators in population1. $n_{2h1}$ and $d$ belong to first phase of
  executing a ATS in $s_{1h}$ and $n_1$ and $d_1$ belong to first phase of executing a ATS in all a
  PSU. m is number of PSU that is selected in the first stage.
 \label{tab}}
 % \end{center}
%\begin{center}
\begin{tabular}{c|ccc|ccc}
eff, $n_{1h}$=50&$c_{tar}/c_{aux}$=10&&&$c_{tar}/c_{aux}$=5&&\\
$m$&4&3&4&4&3&4\\ \hline
$(n_{2h1},d,n_{1},d_1)$& (10,4,13,10)&(10,4,13,10)&(6,4,9,10)&(10,4,16,10)&(10,4,16,10)&(5,4,12,10)\\

\hline
$RegO_x$   & 2.697 &  1.724 &  2.757 &  2.021 &  1.405  & 1.793\\
$Reg1_x$   & 2.610  &  1.692 &  2.660 &   1.956&   1.376 &  1.730\\
$Regopt_x$ & 2.260  & 1.654  &2.145  &1.668  & 1.351  & 1.320\\
$Regb1_x$   &2.208 &  1.652 &  2.082 &  1.628  & 1.343   &1.279\\
$ATSC_x$   &1.087  & 1.010  &  1.070   & 1.008 &  1.000  & 1.086\\
$Regs_x$    &1.766  &1.557  & 1.588  & 1.258  & 1.254  & 1.010\\
\hline
%$n_{1h}$=50&$c_{tar}/c_{aux}$=10&&&$c_{tar}/c_{aux}$=5&&\\
$(n_{2h1},d,n_{1},d_1)$& (10,3,13,12)&(10,3,13,12)&(6,3,9,12)&(9,3,16,12)&(9,3,16,12)&(5,3,12,12)\\
%$m$&4&3&4&4&3&4\\
\hline
$RegO_z$   & 1.403  & 1.218  & 1.158  & 1.025  & 0.925 &  0.712\\
$Reg1_z$   & 1.438  & 1.242 &  1.186 & 1.053  & 0.946  & 0.730\\
$Regopt_z$ & 1.051  & 1.151 &  0.961 &  0.740  &  0.827  & 0.523\\
$Regb1_z$  & 1.391  & 1.328  & 1.582 &  1.080   &1.057  & 1.094\\
$ATSC_z$  & 1.236  & 1.159  & 1.182 &  1.372  & 1.186 &  1.266\\
$Regs_z$   & 1.116  & 1.168  &1.372 &  0.850   & 0.903 &  0.987\\
\hline
%$n_{1h}$=50&$c_{tar}/c_{aux}$=10&&&$c_{tar}/c_{aux}$=5&&\\
$(n_{2h1},d,n_{1},d_1)$& (10,4,12,10)&(10,4,12,9)&(6,4,8,10)&(10,4,16,9)&(10,4,16,10)&(5,4,11,10)\\
%$ps$&4&3&4&4&3&4\\
\hline
$RegO_w$  &  1.481   &1.26  &  1.253&   1.071 &  1.017 &  0.741\\
$Reg1_w$    &1.489  & 1.264 &  1.256  & 1.074 &  1.018 &  0.745\\
$Regopt_w$  &1.012 &  1.083  & 1.134&   0.736  & 0.868&   0.732\\
$Regb1_w$  & 1.014   &1.079   &1.137  & 0.738 &  0.865 &  0.734\\
$ATSC_w$   &1.096   &1.015   &1.036 &  1.134&   1.091 &  1.049\\
$Regs_w$ &   0.854 &  0.919   &0.949  &0.630  &  0.753 &  0.648\\
\hline
 rbias&&&&&&\\

$(n_{2h1},d,n_{1},d_1)$& (10,4,13,10)&(10,4,13,10)&(6,4,9,10)&(10,4,16,10)&(10,4,16,10)&(5,4,12,10)\\
\hline
$RegO_{x}$  &  0.002  & 0.016  & -0.007  &-0.002&  -0.001 & 0.000   \\
$Reg1_{x}$   & 0.003  & 0.017  & -0.007 & -0.002 & -0.001 & 0.000   \\
$Regopt_{x}$ & -0.066 & -0.066 & -0.171 & -0.074 & -0.079 & -0.216  \\
$Regb1_{x}$  & -0.067 & -0.073&  -0.174  &-0.076 & -0.085  &-0.219  \\
$ATSC_{x}$  & -0.009 & 0.011  & 0.001  & -0.006 & 0.003  & -0.006  \\
$Regs_{x}$   & -0.110  & -0.119 & -0.209 & -0.119 & -0.127 & -0.285  \\
 \hline
$(n_{2h1},d,n_{1},d_1)$& (10,3,13,12)&(10,3,13,12)&(6,3,9,12)&(9,3,16,12)&(9,3,16,12)&(5,3,12,12)\\

\hline

$RegO_{z}$  &  0.024  & -0.003&  -0.007 & -0.003 & -0.008 & -0.001  \\
$Reg1_{z}$  &  0.022  &-0.002 & -0.008 &-0.003  &-0.007 & -0.001  \\
$Regopt_{z}$ &-0.038 & -0.052 & -0.111 & -0.066  &-0.057  &-0.115  \\
$Regb1_{z}$  &-0.094 & -0.124 & -0.235 & -0.138  &-0.153 & -0.279  \\
$ATSC_{z}$  & 0.009  & -0.012 & 0.007  & -0.01  & -0.001&  0.001   \\
$Regs_{z}$  &  -0.182 & -0.208 & -0.312 & -0.203 & -0.225 & -0.355  \\

\hline
%$n_{1h}$=50&$c_{tar}/c_{aux}$=10&&&$c_{tar}/c_{aux}$=5&&\\
$(n_{2h1},d,n_{1},d_1)$& (10,4,12,10)&(10,4,12,9)&(6,4,8,10)&(10,4,16,9)&(10,4,16,10)&(5,4,11,10)\\
%$ps$&4&3&4&4&3&4\\
\hline

$RegO_{w}$ &   -0.001&  -0.005&  0.000 &  0.001 &  0.010 &   -0.009  \\
$Reg1_{w}$ &   -0.002 & -0.006 & 0.000  & 0.001  & 0.010 &   -0.008  \\
$Regopt_{w}$ & -0.111 & -0.147 & -0.280   &-0.113  &-0.125&  -0.342  \\
$Regb1_{w}$  & -0.113 & -0.144 & -0.282  &-0.114 & -0.123&  -0.344  \\
$ATSC_{w}$ &  -0.013 & 0.003 &  -0.012  &0.003  & 0.002  & -0.006  \\
$Regs_{w}$   & -0.084&  -0.089 & -0.254  &-0.098 & -0.104  &-0.312  \\

\end{tabular}%
\end{center}

\end{table}

\begin{table}[h]
\begin{center}
  \caption{ efficiency and relative bias of the estimators, Population1.\label{tab}}
  \end{center}
\begin{center}
\begin{tabular}{c|ccc|cc}
eff, $n_{1h}$=70&$c_{tar}/c_{aux}$=10&&&$c_{tar}/c_{aux}$=5&\\
$m$&4&3&4&4&3\\\hline
$(n_{2h1},d,n_{1},d_1)$& (11,5,15,12)&(11,5,15,12)&(14,5,16,13)&(9,4,17,12)&(9,4,17,12)\\

\hline

$RegO_x$   & 3.395   &1.765 &3.317&  2.074 &  1.368\\
$Reg1_x$   & 3.174  & 1.726 &3.082& 1.999  &1.338\\
$Regopt_x$  &2.594  &1.655  &2.755& 1.542  & 1.271\\
$Regb1_x$   &2.454  & 1.646 &2.770& 1.488  &  1.268\\
$ATSC_x$   &1.044  & 0.988 &0.905& 1.028   & 0.975\\
$Regs_x$    &1.839  & 1.541 &1.890&  1.104  & 1.175\\

\hline
$(n_{2h1},d,n_{1},d_1)$& (10,4,15,13)&(10,4,15,13)&(12,4,16,14)&(8,3,17,13)&(8,3,17,13)\\

\hline
$RegO_z$   & 1.500&   1.254  & 1.625&   0.766 &  0.789\\
$Reg1_z$    &1.550  &  1.287&   1.689 &  0.789&   0.804\\
$Regopt_z$  &1.277 &  1.049 &  1.243  & 0.589&   0.636\\
$Regb1_z$   &1.377 &  1.291 &  1.380 &   0.850&    1.001\\
$ATSC_z$  & 1.363 &  1.168 &  1.302 &  1.319 &  1.207\\
$Regs_z$   & 1.068&   1.115&   0.973   &0.714 &  0.861\\

\hline
$(n_{2h1},d,n_{1},d_1)$& (11,5,13,12)&(11,5,13,11)&(14,4,14,14)&(9,4,16,12)&(9,4,16,12)\\

\hline
$RegO_w$    &   1.559   &   1.241   &   1.570    &   0.839   &   0.896\\
$Reg1_w$    &   1.566   &   1.242   &   1.574   &   0.842   &   0.898\\
$Regopt_w$  &   0.963   &   1.025   &   1.032   &   0.603   &   0.770\\
$Regb1_w$   &   0.964   &   1.022   &   1.034   &   0.604   &   0.767\\
$ATSC_w$   &   0.990    &   0.964   &   0.960    &   1.007   &   1.013\\
$Regs_w$    &   0.790    &   0.893   &   0.822   &   0.509   &   0.661\\
\hline rbias&&&&&\\

$(n_{2h1},d,n_{1},d_1)$& (11,5,15,12)&(11,5,15,12)&(14,5,16,13)&(9,4,17,12)&(9,4,17,12)\\

\hline
$RegO_x$  &  -0.001 & 0.002 &  0.004  & -0.002  &-0.004  \\
$Reg1_x$   & -0.001 &0.003  & 0.003   &-0.003  &-0.005  \\
$Regopt_x$  &-0.060  & -0.070 &  -0.027 & -0.096 & -0.108  \\
$Regb1_x$  & -0.065 & -0.078&  -0.032 & -0.100    &-0.119  \\
$ATS_x$  & 0.006  & 0.004  & 0.005  & 0.003  & -0.007  \\
$Regs_x$   &-0.106  &-0.120  & -0.080   &-0.157 & -0.154  \\

\hline
$(n_{2h1},d,n_{1},d_1)$& (10,4,15,13)&(10,4,15,13)&(12,4,16,14)&(8,3,17,13)&(8,3,17,13)\\

\hline

$RegO_z$  & 0.007 &  0.014 &  0.001 &  0.006  & 0.002   \\
$Reg1_z$  & 0.008  & 0.015  & 0.001 &  0.006  & 0.005   \\
$Regopt_z$  &-0.038 & -0.017 & -0.028 & -0.074&  -0.061  \\
$Regb1_z$  & -0.098 & -0.105 & -0.075 & -0.17  & -0.197  \\
$ATS_z$  & -0.005 & 0.012 &  -0.001&  0.008 &  -0.005  \\
$Regs_z$   & -0.210  & -0.202 & -0.002 & -0.245&  0.001   \\

\hline
$(n_{2h1},d,n_{1},d_1)$& (11,5,13,12)&(11,5,13,11)&(14,4,14,14)&(9,4,16,12)&(9,4,16,12)\\

\hline

$RegO_w$   & 0.005 &  0.014 &  0.001 &  -0.002 & 0.002   \\
$Reg1_w$   & 0.005  & 0.014 &  0.001 &  -0.001& 0.002   \\
$Regopt_w$ & -0.088 & -0.111 & -0.048& -0.150  & -0.173  \\
$Regb1_w$  &-0.090  & -0.108 & -0.050   &-0.152 & -0.171  \\
$ATS_w$  & 0.003 &  -0.001 & -0.007 & -0.004 & -0.011  \\
$Regs_w$  &  -0.070 &  -0.081 & -0.038 & -0.120   &-0.114  \\

\end{tabular}%
\end{center}

\end{table}

The results of population2 are in table 4. In the case of high
correlation ($x$) and also for $w$ (with enough sample size) with $c_{tar}/c_{aux}=10$, ATSD is the proper design to
investigate the population. But for the other cases the results
shows SRSWOR with $\bar{y}_s$ is more appropriate. It seems if there is weak correlation between target and auxiliary
variable, because the target variable is not rare, it is better to expanse
all the costs on finding and investigating target variable with
SRSWOR.\\
%And in the case of high correlation we can trust ATSD
%with a proper regression
%estimator.\\
For unbiasedness the results can be summarized as follow. The
amount of unbiasedness is acceptable for almost all the cases.
Also in the cases that ATSD is better than ATS in efficiency,
again bias of $Regopt$ is substantially smaller (or at least
equivalent) than $Regb1$. In high correlation cases, that $ATSD$
is the proper design, $Regb1$ is a little better than $Regopt$ in
efficiency, but as we discussed before, $Regopt$ is better than $Regb1$
according to bias. Then if we look at efficiency and unbiasedness
simultaneously, we prefer to use $Regopt$ is such cases.

\begin{table}[h]
\begin{center}
  \caption{ efficiency and relative bias of the estimators, Population2.\label{tab}}
  \end{center}
\begin{center}
\begin{tabular}{c|ccc|cc}
$n_{1h}$=50&$c_{tar}/c_{aux}$=10&&&$c_{tar}/c_{aux}$=5&\\
$m$&4&3&4&4&3\\\hline
$(n_{2h1},d,n_{1},d_1)$& (7,5,7,11)&(7,5,7,11)&(4,5,5,11)&(7,5,10,10)&(7,5,10,10)\\

\hline

$RegO_x$   & 3.182 &  1.919  & 3.469 &  2.331 &  1.523\\
$Reg1_x$   & 3.211&   1.919&   3.527 &  2.348 &  1.522\\
$Regopt_x$ & 2.691 &  1.605 &  1.893  & 1.900& 1.260\\
$Regb1_x$  & 2.730   & 1.637  & 2.196 &  2.041  & 1.312\\
$ATS_x$  & 0.919   &0.941 &  0.976&   1.056 &  0.975\\
$Regs_x$    &2.220   & 1.499  & 1.870   & 1.653  & 1.217\\

\hline
$(n_{2h1},d,n_{1},d_1)$& (7,5,7,10)&(7,5,7,10)&(4,5,5,10)&(7,5,9,10)&(7,5,9,10)\\

\hline
$RegO_z$   & 1.241 &  1.103  & 0.983  & 0.918 &  0.932\\
$Reg1_z$   & 1.142   &1.053  & 0.927  & 0.846 &  0.890\\
$Regopt_z$ & 0.666  & 0.549  & 0.289  & 0.512&   0.510\\
$Regb1_z$  & 0.949  & 0.848&   0.607 &  0.728 &  0.729\\
$ATS_z$  & 0.744  & 0.823 &  0.797 &  0.746 &  0.845\\
$Regs_z$   & 0.601 &  0.600   &0.322&   0.461  & 0.605\\

\hline
$(n_{2h1},d,n_{1},d_1)$& (7,5,7,10)&(7,5,7,10)&(4,5,5,10)&(7,5,9,10)&(7,5,9,10)\\

\hline
$RegO_w$   &  1.428  & 1.269 &  1.126 &  1.008  & 0.982\\
$Reg1_w$    & 1.447  &1.280  & 1.138  & 1.021  & 0.994\\
$Regopt_w$ & 1.126  & 1.021  & 0.810  & 0.817 &  0.770\\
$Regb1_w$  & 1.135  & 1.020  &  0.839  & 0.822 &  0.771\\
$ATS_w$  & 0.959  & 0.969 &  1.034 &  0.985 &  0.950\\
$Regs_w$   & 0.783  & 0.832 &  0.540 &  0.569 &  0.634\\

\hline rbias&&&&&\\
$(n_{2h1},d,n_{1},d_1)$& (7,5,7,11)&(7,5,7,11)&(4,5,5,11)&(7,5,10,10)&(7,5,10,10)\\

\hline

$RegO_x$   & 0.002 & 0.000   &-0.003 & 0.000  & -0.001  \\
$Reg1_x$    &0.002  & 0.000 & -0.003&  0.000  & -0.001  \\
$Regopt_x$  &-0.004 & -0.02 &  -0.038 & -0.003 & -0.019  \\
$Regb1_x$  & -0.015 & -0.035 & -0.063  &-0.015 & -0.035  \\
$ATS_x$  & 0.000  & -0.005  &-0.007 & 0.002 &  0.002   \\
$Regs_x$  &  -0.022 & -0.048  &-0.047  &-0.020 &  -0.041  \\

\hline
$(n_{2h1},d,n_{1},d_1)$& (7,5,7,10)&(7,5,7,10)&(4,5,5,10)&(7,5,9,10)&(7,5,9,10)\\

\hline

$RegO_z$   &  0.001&   -0.011 & 0.001&   -0.004 &-0.002  \\
$Reg1_z$  &  0.001  & -0.011 & 0.003  & -0.005  &-0.001  \\
$Regopt_z$ & 0.055  & 0.061 &  0.143 &  0.055 &  0.066   \\
$Regb1_z$  & 0.018  & 0.000 &  0.036  & 0.013   &0.008   \\
$ATS_z$   &0.003  & -0.001 & -0.002  &0.004   &0.006   \\
$Regs_z$    &0.034  & 0.024 &  0.072   &0.034 &  0.030    \\

\hline
$(n_{2h1},d,n_{1},d_1)$& (7,5,7,10)&(7,5,7,10)&(4,5,5,10)&(7,5,9,10)&(7,5,9,10)\\

\hline

$RegO_w$  &  0.001 &  -0.001 & 0.006&   0.002&   -0.005  \\
$Reg1_w$    &0.001  & -0.001 & 0.008 &  0.002   &-0.005  \\
$Regopt_w$ & -0.008 & -0.032  &-0.055&  -0.007  &-0.035  \\
$Regb1_w$  & -0.007 &-0.033  &-0.059&  -0.006 & -0.035  \\
$ATS_w$  &-0.004 & 0.006   &0.004 &  0.009  & 0.009   \\
$Regs_w$   & -0.002  &-0.026 & -0.006 & 0.004 &  -0.015  \\

\end{tabular}%
\end{center}

\end{table}

\section{Conclusion}
ATSD is double sampling version of ATS that can be useful to
investigate rare and cluster population with presenting auxiliary
variables. The results in the simulations are conditional on the
data sets that we used but they should apply to any population
with similar features. In the case of high correlation the
proposed design has good performance and for middle amount of
correlations it is depend on structure of target variable and
relative costs of target and auxiliary variables. Simulations show
when the variables are rare and relative costs is reasonably high,
the proper strategy is ATSD.
\section{Appendix}
\subsection{Appendix A} We have
\begin{eqnarray*}
E_3(\hat{\mu}_{reg})=E_3(\bar{y}_{n_{2}})+\beta(\bar{x}_{n_{1}}-E_3(\bar{x}_{n_{2}}))=\bar{y}_{n_{1}}+\beta(\bar{x}_{n_{1}}-\bar{x}_{n_{1}})=\bar{y}_{n_{1}}
\end{eqnarray*}
and
\begin{eqnarray*}
E_2E_3(\hat{\mu}_{reg})=E_2(\bar{y}_{n_{1}})=\frac{M}{N}\frac{1}{m}\sum_{h
\epsilon s}t_{{y}_{N_h}}
\end{eqnarray*}
and then

\begin{eqnarray*}
E_1E_2E_3(\hat{\mu}_{reg})=\frac{M}{N}E_1(\frac{1}{m}\sum_{h
\epsilon s}t_{{y}_{N_h}})=\bar{Y}_{N}
\end{eqnarray*}

 Also for part1 we have
\begin{eqnarray*}
part1=V_1E_2E_3(\hat{\mu}_{reg})=\frac{1}{N^2}M^2(1-\frac{m}{M})\frac{S^{2}_{ty_N}}{m}.
\end{eqnarray*}
For part2 we have
\begin{eqnarray*}
V_2E_3(\hat{\mu}_{reg})=V_2(\bar{y}_{n_{1}})=\frac{M^2}{m^2N^2}\sum_{h
\epsilon s}N^2_{h}V_2(\bar{y}_{n_{1h}})=\frac{M^2}{m^2N^2}\sum_{h
\epsilon
s}N^2_{h}(1-\frac{n_{1h}}{N_h})\frac{S^2_{y_{N_h}}}{n_{1h}}
\end{eqnarray*}
and then
\begin{eqnarray*}
part2=E_1V_2E_3(\hat{\mu}_{reg})=\frac{M}{mN^2}\sum_{h=1}^MN^2_{h}(1-\frac{n_{1h}}{N_h})\frac{S^2_{y_{N_h}}}{n_{1h}}.
\end{eqnarray*}
For part3 we have
\begin{eqnarray*}
V_3(\hat{\mu}_{reg})=V_3(\bar{y}_{n_{2}})+\beta^2
V_3(\bar{x}_{n_{2}})-2\beta
C_3(\bar{y}_{n_{2}},\bar{x}_{n_{2}})
\end{eqnarray*}
and
\begin{eqnarray*}
V_3(\bar{y}_{n_{2}})=\frac{M^2}{m^2N^2}\sum_{h \epsilon
s}a^2_hV_3(\hat{t}_{yn_{2h}})
\end{eqnarray*}
where $V_3(\hat{t}_{yn_{2h}})$ is variance of Murthy estimator in
ATS in h-th PSU under $s_{1h}$. Then
\begin{eqnarray*}
E_1E_2V_3(\bar{y}_{n_{2}})=\frac{M}{mN^2}\sum_{h=1}^{M}a^2_hE_2V_3(\hat{t}_{yn_{2h}})
\end{eqnarray*}
And then
\begin{eqnarray*}
part3=E_1E_2V_3(\hat{\mu}_{reg})=\frac{M}{m}
\sum_{h=1}^{M}a^2_hE_2(V_3(\hat{t}_{yn_{2}})
+\beta^2V_3(\hat{t}_{xn_{2}})-2C_3(\hat{t}_{yn_{2}},\hat{t}_{xn_{2}}))
\end{eqnarray*}

\subsection{Appendix B}
For $E(\hat{S}^{2}_{ty_N})$ we have

\begin{eqnarray*}
\hat{S}^{2}_{ty_N}=\frac{1}{2m(m-1)}\sum_{h \epsilon
s}\sum_{h'\neq h}(a_h\hat{t}_{yn_{2h}}-a_{h'}\hat{t}_{yn_{2h'}})^2
\end{eqnarray*}

and with
\begin{eqnarray*}
E_{2,3}(a_h\hat{t}_{yn_{2h}}-a_{h'}\hat{t}_{yn_{2h'}})^2=V_{2,3}(a_h\hat{t}_{yn_{2h}}-a_{h'}\hat{t}_{yn_{2h'}})+
E^2_{2,3}(a_h\hat{t}_{yn_{2h}}-a_{h'}\hat{t}_{yn_{2h'}})\\=V_{2}E_3(a_h\hat{t}_{yn_{2h}}-a_{h'}\hat{t}_{yn_{2h'}})+
E_2V_3(a_h\hat{t}_{yn_{2h}}-a_{h'}\hat{t}_{yn_{2h'}})+(E_{2}E_3(a_h\hat{t}_{yn_{2h}})-E_{2}E_3(a_{h'}\hat{t}_{yn_{2h'}}))^2\\
=V_{2}(a_{h}t_{yn_{1h}}-a_{h'}t_{yn_{1h'}})+E_2(V_3(a_{h}\hat{t}_{yn_{2h}})+V_3(a_{h'}\hat{t}_{yn_{2h'}}))+(E_{2}(a_{h}t_{yn_{1h}})-E_{2}(a_{h'}t_{yn_{1h'}}))^2\\
=(N^2_{h}(1-\frac{n_{1h}}{N_h})\frac{S^2_{y_{N_h}}}{n_{1h}}+N^2_{h'}(1-\frac{n_{1h'}}{N_h'})\frac{S^2_{y_{N_h'}}}{n_{1h'}})
+(a^2_{h}E_2V_3(\hat{t}_{yn_{2h}})+a^2_{h'}E_2V_3(\hat{t}_{yn_{2h'}}))+(t_{y_{N_h}}-t_{y_{N_{h'}}})^2\\
\end{eqnarray*}

we have

\begin{eqnarray*}
E_{2,3}(\hat{S}^{2}_{ty_N})=\frac{1}{m}\sum_{h =
1}^MN^2_{h}(1-\frac{n_{1h}}{N_h})\frac{S^2_{y_{N_h}}}{n_{1h}}I_h+\frac{1}{m}\sum_{h=1}^{M}a^2_hE_2V_3(\hat{t}_{yn_{2h}})I_h\\
+\frac{1}{2m(m-1)}\sum_{h =1}^M\sum_{h'\neq
h}^M(t_{yN_h}-t_{yN_h})^2I_{hh'}
\end{eqnarray*}
then

\begin{eqnarray*}
E(\hat{S}^{2}_{ty_N})=\frac{1}{M}\sum_{h =
1}^MN^2_{h}(1-\frac{n_{1h}}{N_h})\frac{S^2_{y_{N_h}}}{n_{1h}}+\frac{1}{M}\sum_{h=1}^{M}a^2_hE_2V_3(\hat{t}_{yn_{2h}})\\
+\frac{1}{2m(m-1)}\frac{m(m-1)}{M(M-1)}2M\sum_{h
=1}^M(t_{yN_h}-\bar{t})^2
\end{eqnarray*}

and finally we have

\begin{eqnarray*}
E(\frac{1}{N^2}M^2(1-\frac{m}{M})\frac{\hat{S}^{2}_{ty_N}}{m})=\frac{1}{N^2}[M^2(1-\frac{m}{M})\frac{S^{2}_{ty_N}}{m}+\frac{M}{m}\sum_{h=1}^{M}
N^2_{h}(1-\frac{n_{1h}}{N_h})\frac{S^2_{y_{N_h}}}{n_{1h}}\\+\frac{M}{m}\sum_{h=1}^{M}a^2_hE_2V_3(\bar{y}_{n_{2h}})-\sum_{h=1}^{M}
N^2_{h}(1-\frac{n_{1h}}{N_h})\frac{S^2_{y_{N_h}}}{n_{1h}}-\sum_{h=1}^{M}a^2_hE_2V_3(\hat{t}_{yn_{2h}})].
\end{eqnarray*}
Now for estimating $S^2_{x_{N_h}}$ with
\begin{eqnarray*}
\hat{S}^2_{x_{N_h}}=\frac{1}{n_{1h}-1}[\hat{t}_{x^2_{n_{2h}}}-
\frac{\hat{t}^2_{x_{n_{2h}}}}{n_{1h}}]
\end{eqnarray*}

we have (for $r=1,2$)

\begin{eqnarray*}
E(\hat{t}_{x^r_{n_{2h}}})=E_{2,3}(\hat{t}_{x^2_{n_{2h}}})=E_2(\sum_{j=1}^{n_{1h}}x^r_{hj})=
\frac{n_{1h}}{N_h}\sum_{j=1}^{N_{h}}x^r_{hj}
\end{eqnarray*}

also

\begin{eqnarray*}
E(\hat{t}^2_{x_{n_{2h}}})=E_{2,3}(\hat{t}^2_{x_{n_{2h}}})=
V_{2,3}(\hat{t}_{x_{n_{2h}}})+E^2_{2,3}(\hat{t}_{x_{n_{2h}}})=
V_2E_3(\hat{t}_{x_{n_{2h}}})+E_2V_3(\hat{t}_{x_{n_{2h}}})+(\frac{n_{1h}}{N_h}\sum_{j=1}^{N_{h}}x_{hj})^2\\
=n^2_{1h}(1-\frac{n_{1h}}{Nh})\frac{S^2_{yN_h}}{n_{1h}}+E_2V_3(\hat{t}_{x_{n_{2h}}})+(\frac{n_{1h}}{N_h}\sum_{j=1}^{N_{h}}x_{hj})^2
\end{eqnarray*}

therefore we have

\begin{eqnarray*}
E(\hat{S}^2_{x_{N_h}})=S^2_{x_{N_h}}-\frac{1}{n_{1h}(n_{1h}-1)}E_2V_3(\hat{t}_{x_{n_{2h}}})
\end{eqnarray*}

 Now with all above computation, an asymptotic unbiased estimator for variance of
the estimator is (if we set $\beta$ instead of $\hat{\beta}$ the
estimator will be unbiased):
\begin{eqnarray*}
\hat{var}(\hat{\mu}_{reg})=\frac{1}{N^2}[M^2(1-\frac{m}{M})\frac{\hat{S}^{2}_{ty_N}}{m}
+\frac{M}{m}\sum_{h\epsilon s}
N^2_{h}(1-\frac{n_{1h}}{N_h})\frac{\hat{S}^2_{y_{N_h}}}{n_{1h}}\\+\frac{M}{m}\sum_{h
\epsilon
s}a^2_h\frac{n_{1h}(N_h-1)}{N_h(n_{1h}-1)}\hat{V}_3(\hat{t}_{y_{n_{2h}}})+\frac{M^2}{m^2}\sum_{h
\epsilon
s}a^2_h(\hat{\beta}^2\hat{V}_3(\hat{t}_{x_{n_{2h}}})-2\hat{\beta}\hat{C}_3(\hat{t}_{y_{n_{2h}}},\hat{t}_{x_{n_{2h}}})]
\end{eqnarray*}

\begin{small}
{}
\end{small}

\end{document}